\documentstyle[12pt]{article}
\begin {document}
\begin{center}
{\large\bf
 A New Condition for the Uniform

 Convergence of certain Trigonometric Series}

S. P. Zhou\footnote{corresponding author\\Supported in part by
Natural Science Foundation of China under grant number
10471130.\\Key wards and phrases: uniform convergence,
quasimonotone, $O$-regularly varying quasimonotone, bounded
varation} and R. J. Le

\end{center}

          \begin{quote}
          \small \bf ABSTRACT.
\rm The present paper proposes a new  condition to replace both the  ($O$-regularly varying) quasimonotone condition and a certain type of bounded variation condition, and shows the same conclusion for the uniform convergence of certain trigonometric series still holds.
          \end{quote}

\begin{center}
\end{center}

\begin{center}
\large\bf \S 1. Introduction
\end{center}

\vspace{3mm}

\rm The following classical result was established by Chaundy and Jolliffe [1]:

\vspace{3mm}

\bf Theorem CJ. \it Suppose that $\{b_{n}\}$ is a non-increasing real sequence with $\lim\limits_{n\to\infty}b_{n}=0$. Then a necessary and sufficient condition for the uniform convergence of the series
$$ \sum_{n=1}^{\infty}b_{n}\sin nx\hspace{.4in}(1)$$
is $\lim\limits_{n\to\infty}nb_{n}=0$.

\vspace{3mm}

\rm Recently, Leindler [2] considered to generalize the monotonicity condition to a certain type of bounded variation condition and proved the following

\vspace{3mm}

\bf Theorem L. \it Let \mbox{\bf b}$=\{b_{n}\}_{n=1}^{\infty}$ be a nonnegative sequence satisfying $\lim\limits_{n\to\infty}b_{n}=0$ and
$$\sum_{n=m}^{\infty}|b_{n}-b_{n+1}|\leq M(\mbox{\bf b})b_{m}\hspace{.4in}(2)$$
for some constant $M(\mbox{\bf b})$ depending only upon $\mbox{\bf b}$ and $m=1, 2, \cdots$. Then a necessary and sufficient condition either for the uniform convergence of series $(1)$, or for the continuity of its sum function $f(x)$, is that $\lim\limits_{n\to\infty}nb_{n}=0$.

\vspace{3mm}

\rm On the other hand, we recall quasimonotonicity. A real sequence $\{b_{n}\}_{n=0}^{\infty}$
is said to be quasimonotone if, for some $\alpha\geq 0$, the sequence $\{b_{n}/n^{\alpha}\}$ is non-increasing. More generally, at least in the form, people can define $O$-regularly varying quasimonotone condition.

For a sequence $\{c_{n}\}_{n=0}^{\infty}$, let
$$\Delta c_{n} = c_{n}-c_{n+1}.$$
A non-decreasing positive sequence $\{R(n)\}_{n=0}^{\infty}$ is said to be $O$-regularly varying if\footnote{In some papers, this requirement is written as that for some $\lambda>1$, $\limsup\limits_{n\to\infty}R([\lambda n])/R(n)<\infty$. It is just a pure decoration.}
$$\limsup_{n\to\infty}\frac{R(2n)}{R(n)}<\infty.\hspace{.4in}(3)$$

A complex sequence $\{c_{n}\}_{n=0}^{\infty}$ is $O$-regularly varying quasimonotone in complex sense if for some $\theta_{0}\in [0,\pi/2)$ and some $O$-regularly varying sequence $\{R(n)\}$ (condition (3)) the sequence
$$\Delta\frac{c_{n}}{R(n)}\in K(\theta_{0}):=\{z: |\mbox{\rm arg} z|\leq\theta_{0}\},\;\;n=0,1,\cdots.$$
Evidently, if $\{c_{n}\}$ is a real sequence, then the $O$-regularly varying quasimonotone condition becomes
$$\Delta\frac{c_{n}}{R(n)}\geq 0,\;\;n=0,1,\cdots,$$\\
which was used in [6] as a generalization of the regularly varying quasimonotone condition and, in particular, the quasimonotone condition.

The following theorem, which also generalizes the classical result of Chaundy and Jolliffe [1], was proved  in Nurcombe [4] in 1992:

\vspace{3mm}

\bf Theorem N. \it If $\{b_{n}\}$ is positive and quasimonotone, then a necessary and sufficient condition either for the uniform convergence of series $(1)$, or for the continuity of its sum function $f(x)$, is that $\lim\limits_{n\to\infty}nb_{n}=0$.

\vspace{3mm}

\rm Theorem N was further generalized by Xie and Zhou [7] to

\vspace{3mm}

\bf Theorem XZ1. \it Let $\{c_{n}\}_{n=0}^{\infty}$ be a complex $O$-regularly varying quasimonotone sequence. Write
$$f(x) := \sum_{k=-\infty}^{\infty}c_{k}e^{ikx}$$
at any point $x$ where the series converges. Suppose
$$c_{n}+c_{-n} \in K(\theta_{1}),\;\;n=1,2,\cdots$$
for some $\theta_{1}\in [0,\pi/2)$. Then the necessary and sufficient conditions for $f\in C_{2\pi}$ and $\lim\limits_{n\to\infty}\|f-S_{n}(f)\|=0$ are that
$$\lim_{n\to\infty}nc_{n} = 0$$
and
$$\sum_{n=1}^{\infty}|c_{n}+c_{-n}| < \infty.$$

\rm As a special case, one has

\vspace{3mm}

\bf Theorem XZ2. \it If $\{b_{n}\}$ is positive and $O$-regularly varying quasimonotone, then a necessary and sufficient condition either for the uniform convergence of series $(1)$, or for the continuity of its sum function $f(x)$, is that $\lim\limits_{n\to\infty}nb_{n}=0$.

\vspace{3mm}

\rm Since quasimonotonicity  and condition (2) are not comparable (cf. [3, Theorem 1]),  one  probably would consider to generalize condition (2) to include both cases. For example, the condition,
$$\sum_{n=m}^{\infty}\left|\frac{b_{n}}{R(n)}-\frac{b_{n+1}}{R(n+1)}\right|\leq M(\mbox{\bf b})\frac{b_{m}}{R(m)}\hspace{.4in}(2')$$
for a nonnegative null sequence and some  $O$-regularly varying sequence $\{R(n)\}$, could be a right one if we still can keep the conclusion of Theorem L. However, after deliberate  investigation, we surprisingly see that a condition $(2^{*})$ (see the statement of Theorem 1 in the next section), a  revision to condition (2) suggested by Leindler, is the right generalization. Furthermore, we can prove that, no matter the general form it looks like, any condition like $(2')$ does imply this new condition $(2^{*})$ (see Section 3)

Throughout the paper, $M(\mbox{\bf c})$ denotes a positive constant depending only upon $\mbox{\bf c}$ (which is independent of $n$ and $x\in [0,2\pi]$) not necessarily the same at each occurrence.
\vspace{3mm}

\begin{center}
\large\bf \S 2. Main Result
\end{center}

\rm Let $C_{2\pi}$ be the space of
 all complex valued continuous functions $f(x)$ of period $2\pi$ with  norm
$$\|f\| = \max_{-\infty<x<+\infty}|f(x)|.$$
Given a trigonometric series $\sum_{k=-\infty}^{\infty}c_{k}e^{ikx}:=\lim\limits_{n\to\infty}\sum_{k=-n}^{n}c_{k}e^{ikx}$, write
$$f(x)=\sum_{k=\infty}^{\infty}c_{k}e^{ikx}$$
at any point $x$ where the series converges. Denote its $n$th partial sum $S_{n}(f,x)$ by
$\sum_{k=-n}^{n}c_{k}e^{ikx}.$

\vspace{3mm}

 \bf Theorem 1. \it Let $\mbox{\bf c}=\{c_{n}\}_{n=0}^{\infty}$ be a complex sequence  satisfying
$$c_{n}\pm c_{-n} \in K(\theta_{0}),\;\;n=1,2,\cdots \hspace{.4in}(4)$$
for some $\theta_{0}\in [0,\pi/2)$. If there is a natural number $N_{0}$ such that
$$\sum_{n=m}^{2m}|\Delta c_{n}|\leq M(\mbox{\bf c})\max_{m\leq n<m+N_{0}}|c_{n}|\hspace{.4in}(2^{*})$$
holds for all $m=1, 2, \cdots$, then the necessary and sufficient conditions for $f\in C_{2\pi}$ and $\lim\limits_{n\to\infty}\|f-S_{n}(f)\|=0$ are that
$$\lim_{n\to\infty}nc_{n} = 0\hspace{.4in}(5)$$
and
$$\sum_{n=1}^{\infty}|c_{n}+c_{-n}| < \infty.\hspace{.4in}(6)$$

\bf Lemma 1 (Xie and Zhou [7, Lemma 2]). \it Let $\{c_{n}\}$ satisfy
$$c_{n}+c_{-n} \in K(\theta_{0}),\;\;n=1,2,\cdots$$
for some $\theta_{0}\in [0,\pi/2)$.
Then $f\in C_{2\pi}$ implies that
$$\sum_{n=1}^{\infty}|c_{n}+c_{-n}| < \infty.$$

\bf Lemma 2. \it Let $\{c_{n}\}$ satisfy all conditions of Theorem $1$. Then $\lim\limits_{n\to\infty}\|f-S_{n}(f)\|=0$ implies that
$$\lim_{n\to\infty}nc_{n} = 0.$$

\bf Proof. \rm Fix $n$, assume
$$\max_{n+jN_{0}\leq k<n+(j+1)N_{0}}|c_{k}|=|c_{k_{j}}|,\;\;j=0, 1, \cdots, [n/N_{0}]-1,$$
$$\hspace{1.0in}n\leq n+jN_{0}\leq k_{j}< n+(j+1)N_{0}\leq 2n.$$
Write
$$S_{4n}(f,x)-S_{n}(f,x) =  \sum_{k=n+1}^{4n}\left(c_{k}e^{ikx}+c_{-k}e^{-ikx}\right)$$
$$= \sum_{k=n+1}^{4n}c_{k}(e^{ikx}-e^{-ikx})
 +\sum_{k=n+1}^{4n}\left(c_{k}+c_{-k}\right)e^{-ikx}.$$
 Let $x_{0}=\frac{\pi}{8n}$, we see that
$$M\sum_{k=n+1}^{4n}\mbox{\rm Re}c_{k}
\leq 2\sum_{k=n+1}^{4n}\mbox{\rm Re}c_{k}\sin kx_{0}
= \left|\sum_{k=n+1}^{4n}\mbox{\rm Re}c_{k}(e^{ikx_{0}}-e^{-ikx_{0}})\right|$$
$$\leq \left|\sum_{k=n+1}^{4n}c_{k}(e^{ikx_{0}}-e^{-ikx_{0}})\right|
\leq \|S_{4n}(f)-S_{n}(f)\|+\sum_{k=n+1}^{4n}\left|c_{k}+c_{-k}\right|.\hspace{.4in}(7)$$
At the same time, it is straightforward to see that $a+b\in K(\theta_{0})$ if $a\in K(\theta_{0})$ and $b\in K(\theta_{0})$ for some $\theta_{0}\in [0,\pi/2]$. Hence we note that condition (4) implies that $|c_{n}|\leq M(\theta_{0})\mbox{\rm Re}c_{n}$ for $n\geq 1$. From condition $(2^*)$, we get for $0\leq j\leq [n/N_{0}]-1$,
$$|c_{2n}|=\left|\sum_{j=2n}^{2(n+jN_{0})-1}\Delta c_{j}+c_{2(n+jN_{0})}\right|\leq \sum_{j=n+jN_{0}}^{2(n+jN_{0})-1}|\Delta c_{j}|+|c_{2(n+jN_{0})}|$$
$$\leq M(\mbox{\bf c})(|c_{k_{j}}|+|c_{2(n+jN_{0})}|)\leq M(\mbox{\bf c}, \theta_{0})(\mbox{\rm Re}c_{k_{j}}+\mbox{\rm Re}c_{2(n+jN_{0})}),$$
therefore,
$$\sum_{k=n+1}^{4n}\mbox{\rm Re}c_{k}\geq M^{-1}(\mbox{\bf c}, \theta_{0})\frac{n}{N_{0}}|c_{2n}|.$$
Applying Lemma 1, with (7), we get $\lim\limits_{n\to\infty}nc_{2n}=0$. The other case $\lim\limits_{n\to\infty}nc_{2n+1}=0$ can be treated similarly. Lemma 2 is proved.

\vspace{3mm}

 \rm The following lemma is a mild variant of a theorem in Paley [5] with almost the same proof.

\vspace{3mm}

\bf Lemma 3. \it Let $f(x)=\sum_{n=1}^{\infty}b_{n}\sin nx \in C_{2\pi}$, and $\{b_{n}\}$ be a complex sequence with $b_{n}\in K(\theta_{0})$ for some $\theta_{0}\in [0,\pi/2)$ and all $n=1,2,\cdots$. Then
$$\lim_{n\to\infty}\|f-S_{n}(f)\|=0.$$

\bf Proof of Theorem 1.

\vspace{2mm}

\it Necessity. \rm Applying Lemmas 1 and 2 we have (5) and (6).

\vspace{2mm}

\it Sufficiency. \rm It is not difficult to see that (with a similar argument to that of (8)) under conditions (5), (6), $\{S_{n}(f,x)\}$ is a Cauchy sequence for each $x$, consequently it converges at each $x$. Now we need only to show that
$$\lim_{n\to\infty}\left\|\sum_{k=n}^{\infty}\left(c_{k}e^{ikx}+c_{-k}e^{-ikx}\right)\right\| = 0\hspace{.4in}(8)$$
in this case. In view of conditions (5) and (6), for given $\epsilon>0$ we choose an $n_{0}>0$ such that for all $n\geq n_{0}$,
$$\sum_{k=n}^{\infty}\left|c_{k}+c_{-k}\right| < \epsilon\hspace{.4in}(9)$$
and
$$n\left|c_{n}\right| < \epsilon.\hspace{.4in}(10)$$
Let $n\geq n_{0}$. Set
$$\sum_{k=n}^{\infty}\left(c_{k}e^{ikx}+c_{-k}e^{-ikx}\right)
=\sum_{k=n}^{\infty}\left(c_{k}+ c_{-k}\right)e^{-ikx}
+2i\sum_{k=n}^{\infty}c_{k}\sin kx$$
$$=: I_{1}(x) + 2iI_{2}(x).$$
By (9),
$$|I_{1}(x)| < \epsilon.$$
Noting for $x=0$ and $x=\pi$ that
$$I_{2}(x) = 0,$$
we may restrict $x$ within $(0,\pi)$ without loss of generality. Take $N=[1/x]$ and set\footnote{When $N\leq n$, the same argument as in estimating $J_{2}$ can be applied to deal with $I_{2}=\sum_{k=n}^{\infty}c_{k}\sin kx$ directly.}
$$I_{2}(x) = \sum_{k=n}^{N-1}c_{k}\sin kx
+\sum_{k=N}^{\infty}c_{k}\sin kx =: J_{1}(x) + J_{2}(x).$$
Now
$$|J_{1}(x)| \leq
 \sum_{k=n}^{N-1}|c_{k}||\sin kx|\leq x\sum_{k=n}^{N-1}k|c_{k}|< x(N-1)\epsilon \leq \epsilon$$
 follows immediately from (10) and $N=[1/x]$. By the following well-known estimate  $$|D_n(x)|=\left|\sum\limits_{k=1}^n\sin kx\right|\leq\frac{\pi}{x},$$
and by Abel's transformation and condition $(2^*)$ (together with (10)),
  $$|J_2(x)|=\left|\sum\limits_{k=N}^\infty c_{k}\sin kx\right|
\leq\sum\limits_{k=N}^\infty|\Delta c_{k}|D_k(x)|+c_{N}|D_{N-1}(x)|$$
$$\leq
Mx^{-1}\left(\sum\limits_{k=N}^\infty|\Delta c_{k}|+|c_{N}|\right).
  \hspace{.4in}(11)$$
  Assume that there is a $k_{j}$, $2^{j}N\leq k_{j}<2^{j}N+N_{0}$, such that  $\max\limits_{2^{j}N\leq k<2^{j}N+N_{0}}|c_{k}|=|c_{k_{j}}|$. We calculate that
$$\sum\limits_{k=N}^\infty|\Delta c_{k}|=\sum\limits_{j=0}^\infty\sum\limits_{2^jN\leq
  k<2^{j+1}N}|\Delta c_{k}|\leq M(\mbox{\bf c})\sum\limits_{j=0}^\infty |c_{k_{j}}|$$
$$= M(\mbox{\bf c})\sum\limits_{j=0}^\infty
k_j |c_{k_j}|k_{j}^{-1}\leq M(\mbox{\bf c})\varepsilon
  N^{-1}\sum\limits_{j=0}^\infty 2^{-j}
\leq
  M(\mbox{\bf c})\varepsilon N^{-1}.  \hspace{.4in}  (12)$$
Altogether, combining (11) and (12) yields that $|J_2(x)|\leq M(\mbox{\bf c})\varepsilon$, and therefore it follows that
$\lim\limits_{n\to
\infty}\|f-S_n(f)\|=0$.

\vspace{3mm}

 From Theorem 1, we have a corollary as follows.

\vspace{3mm}

\bf Theorem 2. \it Let \mbox{\bf b}$=\{b_{n}\}_{n=1}^{\infty}$ be a nonnegative sequence. If there is a natural number $N_{0}$ such that
  $$\sum_{n=m}^{2m}|b_{n}-b_{n+1}|\leq M(\mbox{\bf b})\max_{m\leq n<m+N_{0}}b_{m},\;\; m=1, 2, \cdots,$$
then a necessary and sufficient condition either for the uniform convergence of series $(1)$, or for the continuity of its sum function $f(x)$, is that $\lim\limits_{n\to\infty}nb_{n}=0$.

\bf Proof.  \rm By considering Theorem 1, we only need to show that if $f\in C_{2\pi}$, then
$$\lim_{n\to\infty}\|f-S_{n}(f)\| = 0$$
in this case, and we are done just by applying Lemma 3.

\vspace{3mm}

\begin{center}
\large\bf \S 3. Remark: Any quasimonotonicity implies condition $(2^*)$!
\end{center}

\vspace{3mm}

\rm Obviously, any lacunary trigonometric series satisfies neither condition (2) nor any $O$-regularly varying quasimonotone condition. Theorem 2 is surely a nontrivial generalization of Theorem L. A question is whether condition $(2^*)$ implies some type of quasimonotonicity or on the opposite? After deliberate  investigation, we are surprising to find,  no matter the general form it looks like, any condition even like $(2')$ (it is surely weaker than any quasimonotonicity!) does imply our new condition $(2^{*})$! Therefore Theorem 2 is also an essential  generalization of Theorem XL1 as well as Theorem XL2.

\vspace{3mm}

\bf Theorem 3. \it Let $\mbox{\bf c}=\{c_{n}\}_{n=1}^{\infty}$ be any given complex sequence with $\lim\limits_{n\to\infty}c_{n}=0$ satisfying
$$\sum_{n=m}^{\infty}\left|\frac{c_{n}}{R(n)}-\frac{c_{n+1}}{R(n+1)}\right|\leq M(\mbox{\bf c})\frac{|c_{m}|}{R(m)}\hspace{.4in}(2'')$$
for some $O$-regularly varying sequence $\{R(n)\}$.
 Then $\{c_{n}\}$ satisfies $(2^*)$.

\vspace{3mm}

\bf Proof. \rm Suppose that $\{c_{n}\}$ satisfies $(2'')$. Write
$$\frac{c_{k}}{R(k)}-\frac{c_{k+1}}{R(k+1)}=\frac{c_{k}(R(k+1)-R(k))+R(k)(c_{k}-c_{k+1})}{R(k)R(k+1)},$$
then
$$  c_{k}-c_{k+1}=R(k+1)\left(\frac{c_{k}}{R(k)}-\frac{c_{k+1}}{R(k+1)}\right)-\frac{c_{k}}{R(k)}(R(k+1)-R(k)).$$
On the other hand, for $k\geq n$,
$$\left|\frac{c_{k}}{R(k)}\right|=\left|\sum_{j=k}^{\infty}\left( \frac{c_{j}}{R(j)}-\frac{c_{j+1}}{R(j+1)}\right)\right|
\leq \sum_{j=n}^{\infty}\left|\frac{c_{j}}{R(j)}-\frac{c_{j+1}}{R(j+1)}\right|\leq M(\mbox{\bf c})\frac{|c_{n}|}{R(n)}.$$
Thus under the assumption,
$$ \sum_{n\leq k\leq 2n}|c_{k}-c_{k+1}|\leq R(2n+1) \sum_{n\leq k\leq 2n}\left|\frac{c_{k}}{R(k)}-\frac{c_{k+1}}{R(k+1)}\right|$$
$$\hspace{.2in}+\frac{|c_{n}|}{R(n)} \sum_{n\leq k\leq 2n} (R(k+1)-R(k))$$
$$\leq M(\mbox{\bf c})R(2n+1) \frac{|c_{n}|}{R(n)} +\frac{|c_{n}|}{R(n)}(R(2n+1)-R(n))$$
$$\leq M(\mbox{\bf c, R})|c_{n}|\leq M(\mbox{\bf c, R})\max_{n\leq k<n+N_{0}}|c_{k}|$$
for any fixed natural number $N_{0}$, that is what we are required to prove.

\vspace{3mm}

\bf Corollary. \it Let $\mbox{\bf c}=\{c_{n}\}_{n=1}^{\infty}$ be any given complex $O$-regularly varying quasimonotone sequence with $\lim\limits_{n\to\infty}c_{n}=0$. Then $\{c_{n}\}$ satisfies $(2^*)$.

\vspace{3mm}

\rm From $\sum_{k=1}^{\infty}2^{-\alpha k}\sin 2^{k}x=:\sum_{n=1}^{\infty}b_{n}\sin nx$ we can see that, for any $\alpha>0$, $n^{\alpha}b_{n}\not\to 0$, $n\to\infty$, therefore, the condition that $N_{0}$ is a fixed natural number cannot be removed.

\vspace{3mm}
\newpage
\begin{center}
{\Large\bf References}
\end{center}
\begin{enumerate}
\rm
\item T. W. Chaundy and A. E. Jolliffe, \it The uniform convergence of a certain class of trigonometric series, \rm Proc. London Math. Soc.(2) 15(1916), 214-216.
\item L. Leindler, \it On the uniform convergence and boundedness of a certain class of sine series, \rm Anal. Math. 27(2001), 279-285.
\item L. Leindler, \it A new class of numerical sequences and its applications to sine and cosine series, \rm Anal. Math. 28(2002), 279-286.
 \item J. R. Nurcombe, \it On the uniform convergence of sine series with quasimonotone coefficients, \rm J. Math. Anal. Appl. 166(1992), 577-581.
\item R. E. A. C. Paley, \it On Fourier series with positive coefficients, \rm J. London Math. Soc. 7(1932), 205-208.
\item V. B. Stanojevic, \it $L^{1}$-convergence of Fourier series with $O$-regularly varying quasimonotone coefficients, \rm J. Approx. Theory 60(1990), 168-173.
\item T. F. Xie and S. P. Zhou, \it On certain trigonometric series, \rm Analysis 14(1994), 227-237.

\end{enumerate}

\begin{flushleft}
\rm S. P. Zhou:\\
Institute of Mathematics\\
Zhejiang  Sci-Tech University\\
Xiasha Economic Development Area\\
Hangzhou, Zhejiang 310018  China
\end{flushleft}

\begin{flushleft}
\rm R. J. Le:\\
 Department of Mathematics\\
Ningbo University\\
Ningbo, Zhejiang 315211 China
\end{flushleft}

\begin{center}
\bf Keywords \rm uniform convergence,  quasimonotone,  $O$-regularly varying quasimonotone, bounded variation
\end{center}
\end{document}